\documentclass[12pt]{amsart}

\title{A remark concerning Feynman Kac formulas for the perturbed
 harmonic oscillator }
\author{Lisette Jager}
\email{lisette.jager@univ-reims.fr}

\date{}

\usepackage{graphicx}

\usepackage{amsmath}
\usepackage{amssymb}
\usepackage{amsthm}
 
\newtheorem{theo}{Theorem}
\newtheorem{prop}[theo]{Proposition}
\newtheorem{coro}[theo]{Corollary}
\newtheorem{lem}[theo]{Lemma}

\theoremstyle{remark}

\def\ddd{{\mathcal D}}
\def\hhh{{\mathcal H}}

\def\rrr{{\mathcal R}}

\def\N{{\mathbb N}}
\def\R{{\mathbb R}}

\def\ph{\varphi}

\def\ch{{\rm ch}}
\def\sh{{\rm sh}}

\def\gotm{{\mathfrak M}}

\def\R{{\mathbb R}}
\def\ph{\varphi}

\begin{document}

 \maketitle

\begin{center}{\small
Laboratoire de Math\'ematiques\\ Universit\'e de Reims
Champagne-Ardenne\\ Moulin de la Housse, B. P. 1039\\ F-51687
Reims, France }
\end{center}

\begin{abstract}
We give the solution of certain parabolic evolution problems
 (time-depending perturbations of the heat equation for the 
harmonic oscillator ) as explicit  integrals
on a space of continuous functions, called the Wiener space.  The methods
 are based on the Mehler formula giving the solution of the unperturbed
 problem and on the use of discretization to split the difficulties.
\end{abstract}

\section{Introduction.}

The aim of this article is to give an {\it explicit }
expression
 for the solutions of the  problem
\begin{equation}
\label{pb}
\left\{
\begin{array}{lll}
\displaystyle
\frac{\partial v}{\partial t }(t,x)  -
\frac{\partial^2 v}{\partial x^2 }(t,x)  +(x^2  + c(t,x))v(t,x)
=0 \quad {\rm on }\quad ]0,\infty[\times \R 
 \\
\\ \displaystyle
\lim_{t\rightarrow 0}v(t,x)=v_0(x) 
\end{array}
\right. \ ,
\end{equation}
as a Feynman Kac type integral on the Wiener space $C_W$, which is the space of 
continuous functions on $[0,1]$, equipped with a special measure $m_W$
 called Wiener measure (see Subsection \ref{rappels-Wiener}).
\\
Recall that the theory of semigroups and their perturbations gives
the existence and uniqueness of the solution of (\ref{pb}) in
$C^0(\R^+,L^2(\R))$, if, for example, $c$ is continuous, bounded and satisfies 
the following H\"older condition
\begin{equation} \label{Hoelder-uniqueness}
\exists L>0,\exists \alpha \in [0,1] \ : \ 
\forall s, t\in [0,\infty),\ \forall x\in \R,\
 \ |c(t,x)-c(s,x)|\leq L|t-s|^{\alpha}\ .
\end{equation}
The  problem of giving such explicit solutions has already been
 studied under  more restrictive
 conditions  than in the present case.
 In \cite{Donsker-Lions}, \cite{Kuo}, the case when the potential 
$V(t,x)= x^2  + c(t,x)$
is replaced by a  general, unbounded but time-independant potential $V(x)$ 
is studied and in \cite{Kahane}, the author deals with the case when $V(t,x)$ 
is bounded.\\
The present work is devoted to the study of the mixed case (\ref{pb}).
The main result of this paper   is the following
\begin{theo}\label{solution}
Let $v_0\in L^2(\R)$.
Let $c$ be continuous and bounded on $[0,\infty)\times \R$. Suppose
that $c$ belongs to
$ L^2(]0,T[ \times \R)$ for any $T>0$
and 
 satisfies the H\"older condition (\ref{Hoelder-uniqueness}).
The function $v$ defined on $]0,\infty)\times \R$ by
\begin{equation}\label{def-de-v}
\begin{array}{lll}
v(t,x)= & \displaystyle \int_{C_W} v_0(\sqrt{2t }\ w(1) +x) \
 \exp\left(-t\int_0^1 \left(x +\sqrt{2t}\ w(s)\right)^2 \ ds\right)\\ \\
 & \displaystyle \qquad \qquad \exp\left(-t\int_0^1
c\left(t(1-s),\sqrt{2t}\ w(s)+ x\right)\ ds \right)\ dm_W(w)
\end{array}
\end{equation}
is the  solution of (\ref{pb}) in $C^0(\R^+, L^2(\R))$.
\end{theo}
 This result differs from the usual Feynman-Kac formula, for the 
function $w\in C_W$ is never evaluated in $t$.  No well known
 change of measure (\cite{Kuo}, \cite{Aida}, \cite{Ustunel-Zakai})   on
$C_W$ allows to get one formula from the other one.\\ 
One interest would be to get regularity results (or to study the dependency 
on a parameter) more easily than using the theory of semi-groups, by
straightforward derivations.\\
 Problem  (\ref{pb}) will be considered as a perturbation of
 the heat equation related to the harmonic oscillator $H$.
Section \ref{rappels} contains the notions needed about Wiener integrals,
 followed by some facts about the heat kernel for the harmonic oscillator,
in particular Mehler's formula. The proof of a preliminary version of Theorem
  \ref{solution}, under restrictive regularity conditions, is
 completed in Sections  \ref{suite} and  \ref{soln-pb}
by constructing 
a sequence of functions expressed as 
integrals on $C_W$ and  converging to the solution of (\ref{pb}).
The regularity conditions are weakened in Section \ref{affaiblissement},
 which leads to Theorem \ref{solution}.
In section \ref{another}, a  slight modification of section \ref{suite} gives an alternate demonstration
 of Feynman-Kac formula, which does not rely on It\^ o  calculus.
Eventually  we give some explicit expressions
of Wiener integrals, which can be 
 deduced from Theorem \ref{solution} ( annexe  \ref{explicit-expressions}).


\section{Preliminaries}\label{rappels}

\subsection{Wiener integrals}\label{rappels-Wiener}

The construction of the Wiener measure is detailed, for example, in \cite{Kuo} 
and \cite{Yeh}. The (classical)  {\it Wiener space}, denoted by $C_W$,
 is the set of all
real-valued continuous functions $w$ on $[0,1]$, with $w(0)=0$. It can be
 equipped with a probability measure $m_W$ called {\it Wiener measure},
which is defined on  a $\sigma$-field  $\gotm^*$ containing all
sets of the type
$$
J=\{ w\in C_W \ : \ (w(t_1), \dots, w(t_n))\in \hhh \},
$$
where $ t_0=0< t_1< \dots< t_n \leq 1$ and $\hhh$ is  a Borel set of $\R^n$.
For this kind of set the measure is given by
$$
m_W(J)=    \int_{\hhh} f_{t_1,\dots,t_n}(\xi_1,\dots,\xi_n) 
\ d\xi_1 \dots d\xi_n  , 
$$
where $f_{t_1,\dots,t_n}$ is the {\it normal density}
$$
f_{t_1,\dots,t_n}(\xi_1,\dots,\xi_n)=
\left( (2\pi)^n \prod_{i=1}^n(t_i-t_{i-1}) \right)^{-1/2}
\exp\left( -\frac{1}{2} \sum_{i=1}^{n}
\frac{(\xi_i-\xi_{i-1})^2}{t_i-t_{i-1}}\right)
$$
(with $\xi_0=0$). \\
The integral of an integrable, real-valued function $F$ defined on $C_W$ will be
denoted by $E^W[F]$ or  $\int_{C_W}F(w) \ dm_W(w)$.
When $F(w)$ depends on the values of $w$ at finitely many fixed 
points $0\leq t_1<t_2<\dots t_n$ of $[0,1]$,
 $E^W[F]$ can be expressed as follows.
Let $\ph$ be a measurable real-valued function on $\R^n$ and let $\rrr^1$ denote
 the Borel sets on $\R$.
The mapping
$$
\begin{array}{lllll}
F : & (C_W,\gotm^*)  & \rightarrow & (\R,\rrr^1)\\
    & w & \mapsto & \ph(w(t_1),\dots,w(t_n))\\
\end{array}
$$
is measurable and 
\begin{equation}\label{Wiener-nbfinipoints}
 \int_{C_W}F(w) \ dm_W(w) =
 \int_{\R^n} \ph(\xi_1,\dots,\xi_n)
 f_{t_1,\dots,t_n}(\xi_1,\dots,\xi_n) \ d\xi_1 \dots d\xi_n ,
\end{equation}
in the sense that the existence of one side implies the existence of the other
  and the equality.
When $F$ depends on the value of $w$ at infinitely many values of $t$, 
it is useful to recall that, if the topology on $C_W$ is defined by the
 uniform norm, any open set is  $\gotm^*$-measurable and any real-valued, 
continuous  function is  $\gotm^*$-measurable.\\
One important propertie is that the so-called coordinate process, defined
 on $C_W$ by
$$
\forall t \in [0,1], \ \forall w\in C_W,\ 
X_t(w)= w(t)
$$
is a  Brownian motion. A consequence is Fernique's Theorem, which states
 that there exists a positive $\alpha$ such that the integral
$\displaystyle
\int_{C_W}  e^{\alpha ||w||^2}  \ dm_W(w)  
$ converges, $||w||$  being the  uniform norm of $w$ (\cite{Kuo}).


\subsection{Heat equation for the harmonic oscillator}

Classically  the problem
\begin{equation} \label{chaleur-harmonique}
\left\{ \begin{array}{lll}
\displaystyle
\frac{\partial v}{\partial t }(t,x)  -
\frac{\partial^2 v}{\partial x^2 }(t,x)  +x^2 v(t,x)  
=0\qquad  (t,x)\in ]0,\infty[\times \R,\\
v(0)=v_0
\end{array} \right.
\end{equation}
(where  $v_0\in L^2(\R)$) has a unique solution 
$v$ belonging to  $C^0([0,\infty[, L^2(\R))$ (\cite{Pazy}).
More precisely, there exists a semigroup $(U_t)_{t\geq 0}$ of
$L^2(\R)$-contractions such that, for all $t > 0,$
$v(t,\cdot)= U_tv_0$. In the case of the harmonic oscillator this semigroup 
is explicitely given  by Mehler's formula (\cite{BGV})
\begin{eqnarray}
 (U_tv_0)(x) = 
 \int_{\R} \frac{1}{\sqrt{2\pi \sh(2t)}} 
\exp\left( -\frac{1}{2}\left( (x^2+z^2)\frac{\ch(2t)}{\sh(2t)}
-\frac{2xz}{\sh(2t)} \right)  \right)\  v_0(z)
 \ dz.  \label{Utv1}
\end{eqnarray}
One can check directly that $v : (t,x) \mapsto (U_tv_0)(x)$ is 
infinitely derivable on $]0,\infty[\times \R$.
The change of variable  $  \displaystyle z-\frac{x}{\ch(2t)}=y $ gives 
another expression,
\begin{eqnarray}
 v(t,x) = (U_tv_0)(x) = 
\int_{\R} q(t,x,y)\  v_0\left( y+ \frac{x}{\ch(2t)} \right) \ dy\ 
 \label{Utv2}
\end{eqnarray}
where
$$
q(t,x,y)=(2\pi\ \sh(2t))^{-1/2} 
\exp\left(-\frac{1}{2}\frac{\ch(2t)}{\sh(2t)} y^2 
 -\frac{1}{2} \frac{\sh(2t)}{\ch(2t)} x^2 \right).
$$
Both formulae give solutions of (\ref{chaleur-harmonique}) for
initial conditions which do not necessarily belong to $L^2$.
For example, 
if $v_0$ is continuous and vanishes at infinity, $U_tv_0$ is $C^{\infty} $ on
 $]0,\infty[\times \R$ and continuous on $ [0,+\infty)\times \R$. The 
property
$$
\forall t,\tau \geq 0 \ ,\qquad U_t U_{\tau} v_0 = U_{t+\tau} v_0
$$ 
remains valid as well.


\section{Approximation sequence}
\label{suite}

This Section is the first step in the proof of Theorem \ref{solution}.
It is devoted to the construction of a sequence $(v_n^{(t)})_{n\in \N}$ of 
functions designed to
 approximate the solution of (\ref{pb}). For a given positive integer 
$n$ and a fixed  positive $t$, $n$ factors $e^{\frac{t}{n}H}$ coming from
(\ref{Utv1}) alternate with  $n$ factors  $e^{\frac{t}{n}c(t_k,x)}$
(as in Trotter's formula), where the
$t_k$ s are the bounds of a suitable discretization of $[0,t]$.
 The expression thus found can
be written as an integral on the Wiener space which converges to the 
function $v$ in Theorem \ref{solution}. The exponent $(t)$ stresses the fact
 that the sequence depends on the upper bound of the time interval.
\vskip 0.5 cm 

Suppose that the function $c=c(t,x)$ is continuous and bounded on
 $[0,\infty)\times \R$
and that $v_0=v_0(x)$ is   continuous on $\R$ and vanishes at infinity.
Let us split the interval $[0,t]$ into $2n$ subintervals bounded by
the $ \tau_k= \frac{kt}{2n}, k=0,\dots,2 n$ and assign $v^{(t)}_n$ to
satisfy, on each subinterval, an incomplete version of
$$
\frac{\partial v}{\partial t }(t,x)  -
\frac{\partial^2 v}{\partial x^2 }(t,x)  +x^2 v(t,x)  + c(t,x)v(t,x) =0.
$$
The sequence $v^{(t)}_n$ is constructed by induction in  the following way
\begin{itemize}
\item
 on the  even  interval $[\tau_{2k},\tau_{2k+1}]$, $v^{(t)}_n$ is a
 solution of the heat equation for the harmonic oscillator
\begin{equation}\label{eq-paire}
\frac{\partial v}{\partial \tau }(\tau,x)  -
2\left(\frac{\partial^2 v}{\partial x^2 }(\tau,x)- x^2 v(\tau,x)\right)=0
\end{equation}
and satisfies the initial condition
 $\displaystyle \lim_{\tau \rightarrow \tau_{2k}}
 v(\tau,x)= v_n^{(t)}(\tau_{2k},x)$,
where  $v_n^{(t)}(\tau_{2k},\cdot)$ was built in the preceding step ;
\item
on the odd interval $[\tau_{2k+1},\tau_{2k+2}]$, $v^{(t)}_n$ satisfies
 the ordinary differential equation
\begin{equation}\label{eq-impaire}
\frac{\partial v}{\partial \tau }(\tau,x)+2 c(\tau_{2k+2},x)v(\tau,x)=0,
\end{equation}
with the initial condition
 $\displaystyle \lim_{\tau \rightarrow \tau_{2k+1}}
 v(\tau,x)= v_n^{(t)}(\tau_{2k+1},x)$, 
where  $v_n^{(t)}(\tau_{2k+1},\cdot)$ was built in the preceding step.
\end{itemize}
Both equations have of course an explicit solution
$$
\begin{array}{lll}
\displaystyle v_n^{(t)}(\tau,x)= 
\left(U_{2(\tau-\tau_{2k})}v_n^{(t)}(\tau_{2k},\cdot)\right)(x) &{\rm  on} &
[\tau_{2k},\tau_{2k+1}],\\ \\
\displaystyle  v_n^{(t)}(\tau,x)= e^{-2c(\tau_{2k+2},x)(\tau-\tau_{2k+1} )}
v_n^{(t)}(\tau_{2k+1},x)  &{\rm  on} & [\tau_{2k+1},\tau_{2k+2}].\\
\end{array}
$$
The equation (\ref{eq-impaire}) is a constant coefficient linear differential
 equation since  $c$ depends on the fixed time $\tau_{2k+2}$.
The initial conditions ensure the continuity with respect to $\tau$.
  The factor $2$ in 
 (\ref{eq-paire}) and (\ref{eq-impaire}) compensates the fact that each
 equation
 is solved on half of the interval (see Lemma \ref{Lemme-Kahane} below and
\cite{Kahane}, where this discretization was introduced,  to treat 
the usual heat equation.)

\vskip 0.5 cm

We first write $v_n^{(t)}(t,x)$ as iterated integrals.
\begin{prop}\label{integrales-iterees}
For all real $x$ we have 
$$
\begin{array}{llll}
\displaystyle v^{(t)}_n(t,x)  = \displaystyle
\int_{\R^n}   \left( 2 \pi \sh(2t/n)\right)^{-n/2}
v_0\left(\sigma_n + \frac{x}{\ch(2t/n)^n}\right) \\
\qquad   \displaystyle
\exp\left( -\frac{1}{2} \frac{\ch(2t/n)}{\sh(2t/n)}  
             \sum_{j=1}^n \left( \sigma_{n-j+1}   -
     \frac{\sigma_{n-j}}{\ch(2t/n)}
                    \right)^2   \right)  \\
\qquad  \displaystyle
\exp\left( -\frac{1}{2} \frac{\sh(2t/n)}{\ch(2t/n)}  
             \sum_{j=1}^n \left( \sigma_{n-j}+ 
  \frac{x}{\ch(2t/n)^{n-j}}
                    \right)^2   \right) \\ \\
\qquad  \displaystyle
\exp\left(-\frac{t}{n} \sum_{j=1}^n c\left(\frac{jt}{n},\sigma_{n-j}+ 
  \frac{x}{\ch(2t/n)^{n-j}} \right)  \right)
 \quad  d\sigma_1 \dots d\sigma_n \ 
\end{array}
\label{pertube-vn(t,x)}
$$
 with the convention $\sigma_0=0$.
\end{prop}

\noindent{\it Proof.}
 By induction on $k$ we get that, for all 
 $\tau \in [\tau_{2k},\tau_{2k+1}]$, 
\begin{equation}
\label{vn-intermediaire}
\begin{array}{lllll}
\displaystyle
v^{(t)}_n(\tau,x)= \int_{\R^{k+1}} q(2(\tau-\tau_{2k}),x,y_{k+1}) \\
\\
\displaystyle
\prod_{j =1}^{k} 
q   \left(  t/n, 
\sum_{l=j+1}^{k+1}   \frac{y_{l}}{ \ch(2t/n)^{ l-1-j} }
 + \frac{x}{\ch(4(\tau -\tau_{2k}))  \ch(2t/n)^{ k-j} }, y_j
    \right)
\\ \\
\displaystyle
\prod_{j =1}^{k} 
\exp
\left(-\frac{t}{n} c\left(\tau_{2j}, 
\sum_{l=j+1}^{k+1} 
\frac{y_{l}}{ \ch(2t/n)^{ l-1-j} }
 + \frac{x}{\ch(4(\tau -\tau_{2k}))  \ch(2t/n)^{ k-j} }
\right)
 \right)
\\
\\
\displaystyle 
v_0\left(
\sum_{l=1}^{k+1} 
\frac{y_{l}}{  \ch(2t/n)^{ l-1} }
 + \frac{x}{\ch(4(\tau -\tau_{2k}))  \ch(2t/n)^{ k}}
 \right) 
\quad dy_1\dots dy_{k+1}.
\end{array}
\end{equation}
For $\tau=\tau_{2n}=t$ and $k+1=n$, this gives
$$
\begin{array}{lllll}
\displaystyle
v^{(t)}_n(t,x)= \int_{\R^{n}} 
v_0\left(   \sum_{l=1}^{n} 
\frac{y_{l}}{ \ch(2t/n)^{ l-1}} + \frac{x}{ \ch(2t/n)^{n} }\right) 
 \\ \\  \displaystyle
\prod_{j =1}^{n} 
q\left( t/n , 
\sum_{l=j+1}^{n} 
\frac{y_{l}}{\ch(2t/n)^{l-1-j}}+\frac{x}{\ch(2t/n)^{n-j}},y_{j}\right)
\\ \\ \displaystyle
\prod_{j =1}^{n}  \exp \left(
-\frac{t}{n}c\left(\tau_{2j}, \sum_{l=j+1}^{n} 
\frac{y_{l}}{ \ch(2t/n)^{ l-1-j} } + \frac{x}{  \ch(2t/n)^{ n-j} } \right)
 \right)
\quad  dy_1\dots dy_{n},
\end{array}
$$
with the convention that a sum is equal to zero if its lower index  is strictly
superior to its upper index (this is useful for  $j=n$). Let
\begin{equation} \label{chgtvar-sigma}
\sigma_{k}  =
\sum_{l=n+1-k}^{n} \frac{y_{l}}{ \ch(2t/n)^{ l-1-n+k}}\ , \ 1\leq k\leq n
\end{equation}
(recall $\sigma_0=0$). Reciprocally,
$$
y_j=\sigma_{n+1-j}- \frac{\sigma_{n-j}}{\ch(2t/n)}\ ,
j=1,\dots, n.
$$
and this change of variables $y \rightarrow \sigma $
 has Jacobian equal to $1$. Hence
$$
\begin{array}{lllll}
\displaystyle
v^{(t)}_n(t,x)= \int_{\R^{n}} 
v_0\left( \sigma_n
 + \frac{x}{ (\ch(2t/n))^n}
 \right)  \\
\\
\displaystyle
\prod_{j =1}^{n} 
q\left( \frac{t}{n}, 
\sigma_{n-j}
 + \frac{x}{(\ch(2t/n))^{n-j}} , 
\sigma_{n+1-j}- \frac{\sigma_{n-j}}{\ch(2t/n)},
     \right)
\\ \\
\displaystyle
\prod_{j =1}^{n} 
\exp
\left(-\frac{t}{n}c\left(jt/n, 
\sigma_{n-j}
 + \frac{x}{( \ch(2t/n))^{n-j} }\right) \right)
\quad d\sigma_1\dots d\sigma_{n}.
\end{array}
$$
Replacing $q$ by its expression gives the Proposition. \hfill $\square$

\vskip 0.5 cm

The second step is the transformation of this integral on $\R^n$ into an 
integral on the Wiener space. Since $C_W$ does not depend on $n$ we then shall 
be able to let $n$ converge to infinity. These computations lead to the 
following Proposition. 

\begin{prop}\label{expr-V-en-t}
Let $v_0$ be  continuous and bounded on $\R$.
Suppose $c$ is continuous and has a lower bound on  $[0,\infty)\times \R$. 
Moreover let $c$ have an $x$-derivative $c_x$ at each point of
 $]0,\infty)\times\R$ and suppose $c_x$ has an extension to 
$[0,\infty)\times\R$ which is continuous and bounded.\\
When $n$ goes to infinity the  numerical sequence 
$v^{(t)}_n(t,x)$ converges to a limit $v(t,x)$ defined by
$$
\begin{array}{lll}
\displaystyle
v(t,x)=  \int_{C_W} v_0(\sqrt{2t } \ w(1) +x) \ \exp\left( -t 
\int_0^1 \left(  x +  \sqrt{2t}\   w(s) \right)^2 \ ds\right)\\ \\
\displaystyle \qquad \qquad \exp\left(-t\int_0^1
c\left( t(1-s), \sqrt{2t} w(s)+ x\right)  \ ds \right)   \ dm_W(w).
\end{array}
$$
\end{prop}

\noindent{\it Proof.}
To bring out the normal density one uses the change of variables
\begin{equation} \label{chgtvar-xi}
\xi_j=\sigma_j(n \ \sh(2t/n))^{-1/2},\ 1\leq j\leq n 
\end{equation}
and the convention $\xi_0=\sigma_0=0$.
It follows
$$
\begin{array}{llll}
\displaystyle v^{(t)}_n(t,x) & = &\displaystyle \int_{\R^n}(2\pi /n)^{-n/2}
v_0\left((n\ \sh(2t/n))^{1/2}\xi_n +\frac{x}{\ch(2t/n)^n}\right) \\ \\
 &&\displaystyle
\exp\left(-\frac{1}{2}(n\ \ch(2t/n)\sum_{\ell=1}^n\left(\xi_{n-\ell+1} -
     \frac{\xi_{n-\ell}}{\ch(2t/n)} \right)^2   \right)  \\ \\
&&\displaystyle
\exp\left( -\frac{1}{2} \frac{\sh(2t/n)}{\ch(2t/n)}  
             \sum_{\ell=1}^n \left( (n\ \sh(2t/n))^{1/2} \xi_{n-\ell}+ 
  \frac{x}{\ch(2t/n)^{n-\ell}}\right)^2   \right) \\ \\
&&\displaystyle
\exp\left( -\frac{t}{n}\sum_{\ell=1}^n
 c\left(\frac{\ell t}{n}, (n\ \sh(2t/n))^{1/2} \xi_{n-\ell}+ 
 \frac{x}{\ch(2t/n)^{n-\ell}}\right) \right)\quad d\xi_1\dots d\xi_n\ .
\end{array}
$$
We introduce the normal density for equidistant points
 $0\leq t_i=i/n\leq 1$ , 
$$
f_{t_1,\dots,t_n}(\xi_1,\dots,\xi_n)=
 (2\pi/n)^{-n/2} 
\exp\left( -\frac{1}{2} \sum_{i=1}^{n}
n(\xi_i-\xi_{i-1})^2\right), 
$$
to make the transition to Wiener space more natural.
We first obtain
$$
\begin{array}{llll}
\displaystyle v^{(t)}_n(t,x)  =
\displaystyle
\int_{\R^n}  
v_0\left( (n\ \sh(2t/n))^{1/2} \xi_n + \frac{x}{\ch(2t/n)^n}\right)
f_{t_1,\dots,t_n}(\xi_1,\dots,\xi_n) \\
\\
\displaystyle
\exp\left( \frac{n}{2} \sum_{i=1}^{n}
(\xi_i-\xi_{i-1})^2
 -\frac{1}{2}   (n\ \ch(2t/n) 
             \sum_{j=1}^n \left( \xi_{j}   -
     \frac{\xi_{j-1}}{\ch(2t/n)}   \right)^2 
\right) \\ \\
\displaystyle
\exp\left( -\frac{1}{2} \frac{\sh(2t/n)}{\ch(2t/n)}  
             \sum_{j=0}^{n-1} \left( (n\ \sh(2t/n))^{1/2} \xi_{j}+ 
                           \frac{x}{\ch(2t/n)^{j}}   \right)^2    \right) 
\\ \\
\displaystyle
\exp\left( -\frac{t}{n}
             \sum_{j=0}^{n-1}
 c\left(\frac{(n-j) t}{n}, (n\ \sh(2t/n))^{1/2} \xi_{j}+ 
  \frac{x}{\ch(2t/n)^{j}} \right) 
                    \right) 
 \quad    d\xi_1 \dots d\xi_n \ ,
\end{array}
$$
and then, by formula
(\ref{Wiener-nbfinipoints}),
$$
\begin{array}{llll}
\displaystyle v^{(t)}_n(t,x)  =
\int_{C_W}  
F_n(w)G_n(w) H_n(w)\ 
 v_0\left( (n\ \sh(2t/n))^{1/2} w(1) + \frac{x}{\ch(2t/n)^n}\right)
 \quad 
dm_W(w),
\end{array}
$$
with
$$
\begin{array}{lll}
\displaystyle
F_n(w)= 
\exp\left( \frac{n}{2}    \sum_{i=1}^{n}
\left[(1-\ch(2t/n))w(t_i)^2 + (1-\ch(2t/n)^{-1})w(t_{i-1})^2\right]
\right)
\\  \\
\displaystyle
G_n(w)=
\exp\left( -\frac{1}{2} \frac{\sh(2t/n)}{\ch(2t/n)}  
             \sum_{j=0}^{n-1} \left( (n\ \sh(2t/n))^{1/2} w(t_{j}) + 
  \frac{x}{\ch(2t/n)^{j}}
                    \right)^2   \right) 
\\  \\
\displaystyle
H_n(w)=
\exp\left( -\frac{t}{n}
             \sum_{j=0}^{n-1}
 c\left(\frac{(n-j) t}{n}, (n\ \sh(2t/n))^{1/2} w(j/n)+ 
  \frac{x}{\ch(2t/n)^{j}} \right) 
                    \right) .
\end{array}
$$
The Proposition is now a consequence of the dominated convergence Theorem. The
convergences and estimations concerning the four factors are treated in the
Lemma just below
  
\begin{lem}\label{convergences}
Denote by  $m$ a lower bound of $c$.
For all $w\in C_W$,\\
\begin{itemize}
\item
$0\leq F_n(w)\leq 1   \quad  { and} \quad 
\displaystyle \lim_{n\rightarrow \infty} F_n(w)=1 ,$\\
\item
$0\leq G_n(w)\leq 1   \quad  { and} \quad 
\displaystyle
\lim_{n\rightarrow \infty} G_n(w)=
\exp\left( -t\int_0^1 \left(x+\sqrt{2t}\ w(s)\right)^2 \ ds \right),$\\
\item
$0\leq H_n(w)\leq e^{-m t}
 \ \ { and } \ \
\displaystyle
\lim_{n\rightarrow \infty} H_n(w)= 
\exp\big(-t\int_0^1
c\left( t(1-s), \sqrt{2t}\ w(s)+ x\right)  \ ds \big) ,$\\
\item
$\displaystyle
\left| v_0\left((n\ \sh(2t/n))^{1/2} w(1)+\frac{x}{\ch(2t/n)^n}\right)\right|
\leq ||v_0||_{\infty}$\
and \\
$
\displaystyle
\lim_{n\rightarrow \infty} 
 v_0\left( (n\ \sh(2t/n))^{1/2} w(1) + \frac{x}{\ch(2t/n)^n}\right)
= v_0(\sqrt{2t}\  w(1) + x).$\\
\end{itemize}
\end{lem}

\noindent{\it Proof  of Lemma \ref{convergences} }\\
We shall use more than once the continuity of  $w\in C_W$ to compute limits
of Riemann sums like  $\lim_{n\rightarrow \infty}\sum_{j=0}^{n-1} w(t_j)$.
The positivity and the estimates proposed are obvious, except for
$F_n$ which we shall treat first.\\
Let us develop the argument of the exponential
$$
A_n=\frac{n}{2}    \sum_{i=1}^{n}
\left[(1-\ch(2t/n))w(t_i)^2 + (1-\ch(2t/n)^{-1})w(t_{i-1})^2\right].
$$
Since  $t_n=1$ and $w(t_0)=w(0)=0$,
$$
A_n= \frac{n}{2} (1-\ch(2t/n))w(1)^2 -
\frac{n}{2}\frac{ ( \ch(2t/n)-1)^2  }{\ch(2t/n)}
  \sum_{i=1}^{n-1}  w(t_i)^2 \leq 0 \ ,
$$
and then $F_n\leq 1$.\\
Now we claim that $A_n$ converges to  $0$. An asymptotic expansion is enough for
the term containing $w(1)$. The second term can be written as
$$
\frac{n^2}{2}\frac{ ( \ch(2t/n)-1)^2  }{\ch(2t/n)}\  \times \
 \frac{1}{n} \sum_{i=1}^{n-1}  w(t_i)^2.
$$
In this product
the second factor converges to  $\int_0^1 w^2 \ ds$ and the first one
to $0$, using an asymptotic expansion.
\vskip 0.3 cm
To compute the limit of $G_n$ we split the argument of the exponential into 
three terms :
$$
\begin{array}{lll}
\displaystyle
A_n=  -\frac{1}{2} \frac{\sh(2t/n)}{\ch(2t/n)}  
             \sum_{j=0}^{n-1} (n\ \sh(2t/n)) w(t_{j})^2    \\ \\
\displaystyle
B_n=  - \frac{\sh(2t/n)}{\ch(2t/n)}  
             \sum_{j=0}^{n-1}  (n\ \sh(2t/n))^{1/2} w(t_{j})  
  \frac{x}{\ch(2t/n)^{j}}    \\ \\
\displaystyle
C_n=  -\frac{1}{2} \frac{\sh(2t/n)}{\ch(2t/n)}  
             \sum_{j=0}^{n-1} \frac{x^2}{\ch(2t/n)^{2j}} \ . \\ \\
\end{array}
$$
A direct computation gives
$$
A_n= -\frac{n^2}{2}\ \frac{\sh(2t/n)^2}{\ch(2t/n)}  \ \
            \frac{1}{n} \sum_{j=0}^{n-1}  w(t_{j})^2  \quad
\longrightarrow \quad   -\frac{(2t)^2}{2}    \int_0^1 w^2(s) \ ds  \ .
$$
The third term is a geometric sum. An asymptotic expansion leads to
$$
C_n\  \longrightarrow \  -tx^2.
$$ 
The second term is decomposed as 
$$
\begin{array}{lcc} 
\displaystyle
B_n  & = 
\displaystyle
\underbrace{ -x \ \frac{n^{3/2}\sh( \frac{2t}{n})^{3/2}}{\ch( \frac{2t}{n})}  \
            \frac{1}{n} \sum_{j=0}^{n-1}  w(t_{j})} & \displaystyle 
\underbrace{
 -x\  \frac{n^{3/2}\sh( \frac{2t}{n})^{3/2}}{\ch( \frac{2t}{n})}  \
    \frac{1}{n}  \sum_{j=0}^{n-1}  w(t_{j}) (\ch( \frac{2t}{n})^{-j} -1))}\\
 & D_n & E_n
\end{array} .
$$
Since
$$
\frac{n^{3/2}\sh(2t/n)^{3/2}}{\ch(2t/n)} =(2t)^{3/2}(1+o(1/n)),
$$
one has
$$
D_n \ \longrightarrow\ -x (2t)^{3/2}\int_0^1 w(s) \ ds.
$$
Moreover
$$
0\leq 1-  \frac{1}{\ch(2t/n)^{j}}  \leq 1-  \frac{1}{\ch(2t/n)^{n}}
\sim 2 t^2 n^{-1} 
$$
then
$$
|E_n|\leq 
 |x|  \frac{n^{3/2}\sh(2t/n)^{3/2}}{\ch(2t/n)}  \ \ 
 (1-\ch(2t/n)^{-n} )) \ \ 
            \frac{1}{n}  \sum_{j=0}^{n-1} | w(t_{j})|
 \quad \longrightarrow \quad 0.
$$
To sum up,  $(\ch(2t/n))^{-j}$ can be replaced by $1$ in all the
 $w(t_j) (\ch(2t/n))^{-j}$. We conclude that
$$
\begin{array}{ccc}
A_n + B_n + C_n  & \longrightarrow  & \displaystyle
  - 2t^2   \int_0^1 w^2(s) \ ds  \quad
-x (2t)^{3/2}\int_0^1 w(s) \ ds\quad   -tx^2 \\ \\
&=&\displaystyle -t\int_0^1 \left(x+(2t)^{1/2}w(s)\right)^2\ ds. \\
\end{array}
$$

\vskip 0.3 cm
Let us turn to $H_n$.
For all $j$ and all  $w$ we can write
$$
- c\left(\frac{(n-j) t}{n}, (n\ \sh(2t/n))^{1/2} w(j/n)+ 
  \frac{x}{\ch(2t/n)^{j}} \right) \leq -m,
$$
hence the estimate $H_n(w)\leq e^{-mt}.$
Denote by $||c_x||_{\infty}$ the uniform norm of $c_x$ on
$\in[0,\infty)\times\R$. The mean value theorem implies that 
$$
\begin{array}{lll}
\displaystyle
\left| c\left( t-\frac{j}{n}t , (n\ \sh(2t/n))^{1/2} w(j/n)+ 
  \frac{x}{\ch(2t/n)^{j}} \right)-
c\left( t-\frac{j}{n}t, \sqrt{2t}\  w(j/n)+ x\right) \right|
 \\
\quad \displaystyle \leq
||c_x||_{\infty} \left| 
 (n\ \sh(2t/n))^{1/2} w(j/n) - \sqrt{2t}\ w(j/n)   \right|+
||c_x||_{\infty}
 \left|  \frac{x}{\ch(2t/n)^{j}} - x \right|
 \\
\quad \displaystyle \leq
||c_x||_{\infty} ||w||_{\infty} \left| 
 (n\ \sh(2t/n))^{1/2} - \sqrt{2t}    \right|+
||c_x||_{\infty}|x|
 \left|  \frac{1}{\ch(2t/n)^{n}} - 1 \right|,
\end{array}
$$
in which 
the last term is independant of $j$. Therefore
$$
\begin{array}{lll}
\displaystyle
\frac{t}{n}
\left| \sum_{j=0}^{n-1}
 c\big( t-\frac{j}{n}t , (n\ \sh( \frac{2t}{n}))^{1/2} w(j/n)+ 
  \frac{x}{\ch( \frac{2t}{n})^{j}} \big)- 
\sum_{j=0}^{n-1}
c\big( t-\frac{j}{n}t, \sqrt{2t}\ w(j/n)+ x\big) \right|
\\
\qquad \displaystyle \leq
 t ||c_x||_{\infty}
\big(  ||w||_{\infty} \left| 
 (n\ \sh(2t/n))^{1/2} - \sqrt{2t}    \right|+
|x|
 \left|  \frac{1}{\ch(2t/n)^{n}} - 1 \right|\big)\ ,
\end{array}
$$
which converges to $0$ when $n$ goes to infinity.
We deduce that $H_n(w)$ and
$$\displaystyle
\exp\left(-\frac{t}{n}  \sum_{j=0}^{n-1}
c\left( t-\frac{j}{n}t, \sqrt{2t}\  w(j/n)+ x\right)  \right) 
$$
converge to the same limit.
Since the function $u\mapsto c\left( t-ut, \sqrt{2t}\  w(u)+ x\right)$
is continuous on  $[0,1]$, this limit is
$$
\exp\left(-t\int_0^1
c\left( t(1-s), \sqrt{2t}\  w(s)+ x\right)  \ ds \right).
$$
\vskip 0.3 cm
Let us treat the last point of the Lemma. Recall that $v_0$ is continuous
 and bounded. Its argument goes to $(2t)^{1/2} w(1) + x$ because
$$
\ch(2t/n)^n = \exp\left( \frac{2t^2}{n} + o(\frac{1}{n^2}) \right)
\rightarrow 1 ,
$$
which completes the proof. \hfill $\square$


\section{Preliminary version of Theorem \ref{solution} }
\label{soln-pb}

We still need to show that the function $v$ constructed above as the
limit, at time $t$, of the sequence $(v_n^{(t)})_{n\in \N}$ is a solution 
of (\ref{pb}). The demonstration requires stronger regularity conditions on
 $v_0$ and $c$ than the ones used to compute the limit. 
Here is the ( weaker ) version of Theorem \ref{solution} which will be proved in
this section.

\begin{theo}\label{solution-weak}
Let $v_0$ be a  $C^4$ function over $\R$, which has bounded derivatives of 
order up to $4$.
 Suppose $v_0(x)$ converges to $0$ when $x$ goes to infinity.\\
Let $c$ be a function which
\begin{itemize}
\item   is continuous and bounded on  $[0,\infty)\times \R$, 
\item
is  $C^1 $ on  $]0,\infty)\times \R$,
\item
has  bounded space derivatives, up to order $4$, these derivatives being
continuous and bounded on  $]0,\infty)\times \R$.
\end{itemize}
The function $v$ defined on $]0,\infty)\times \R$ by
$$
\begin{array}{lll}
v(t,x)= & \displaystyle \int_{C_W} v_0(\sqrt{2t }\ w(1) +x) \
 \exp\left(-t\int_0^1 \left(x +\sqrt{2t}\ w(s)\right)^2 \ ds\right)\\ \\
 & \displaystyle \qquad \qquad \exp\left(-t\int_0^1
c\left(t(1-s),\sqrt{2t}\ w(s)+ x\right)\ ds \right)\ dm_W(w)
\end{array}
$$
is a solution of (\ref{pb}).
\end{theo}

We  need  preliminary results. The first Lemma shows that the sequence
 $(v_n^{(t)})$  and one of its derivatives 
converge on a dyadic subset of $[0,t]$. The  second one 
gives uniform estimates concerning some of the derivatives
of $v_n^{(t)}$. The  last result proves that a subsequence of  $(v_n^{(t)})$
 converges {\it uniformly}, and  on $[0,t]$ itself.

\begin{lem}\label{cv-sur-D} 
Let $t\in ]0,\infty[$, let  $\ddd$ be the set 
$$
\ddd= \Big\{\tau \in [0,t] \ : \ \exists n_0 \in \N,
 \exists k_0\in \{0,\dots, 2^{n_0}-1\}, \tau= \frac{k_0}{2^{n_0}} t.\Big\}
$$
For all $\tau \in \ddd$ and all $x\in \R$, 
$$
\lim_{n\rightarrow \infty} v_{2^n}^{(t)}(\tau,x) = v(\tau,x) \quad
{\rm and }\quad 
\lim_{n\rightarrow \infty}
\frac{\partial^2 v_{2^n}^{(t)}}{\partial x^2} (\tau,x) =
\frac{\partial^2 v}{\partial x^2}(\tau,x),
$$
where $v$ is the function defined in Theorem \ref{solution}.
\end{lem}

\vskip 0.2 cm
\noindent
The rest of this Section is devoted to the proof of the Lemmas and of
Theorem  \ref{solution-weak}.
\vskip 0.2 cm

\noindent{\it Proof  of Lemma \ref{cv-sur-D}}\\
It is easier to express  $v_n^{(t)}(\tau,x)$ when $\tau = kt/n$ is a bound of
the subdivision. Therefore we are led to consider nested subdivisions. A dyadic
 point $\displaystyle\tau= \frac{k}{2^n}t$, which is already a bound of the
subdivision with $2\cdot 2^n$ intervals, is a bound of all following
 subdivisions.  The set $\ddd$ is the unions of all such points.
Formula  (\ref{vn-intermediaire}) gives $v_n^{(t)}(\tau,x)$ for
 $\displaystyle\tau \in [\frac{2k}{2n}t,\frac{2k+1}{2n}t]$.
Adapting the proof of Proposition (\ref{expr-V-en-t}) we deduce the 
expression of $v_{2^{n+p}}^{(t)}(\tau,x)$ when
$\displaystyle\tau= \frac{k}{2^n}t=  \frac{2^pk}{2^{n+p}}t$ and obtain 
\begin{equation}\label{vn-tau-dyadique}
\begin{array}{lllll}
\displaystyle
v_{2^{n+p}}^{(t)}\left(\tau ,x\right)= \int_{C_W} 
 v_0\left(\sqrt{{2^p k}\ \sh(2t/2^{n+p})}w(1)+\frac{x}{\ch(2t/2^{n+p})^{{2^p k}}}
    \right)
\\ \\
\displaystyle
\exp \left( -\frac{1}{2} \frac{\sh(2t/2^{n+p})}{\ch(2t/2^{n+p})}
 \sum_{l =0}^{{2^p k -1}}
 \left(  \sqrt{{2^p k}\ \sh(2t/2^{n+p})  }   w\left(\frac{l}{{2^p k}}\right)  
 + \frac{x}{\ch(2t/2^{n+p})^{l}} \right)^2
  \right)
\\ \\
\displaystyle
\exp \left( \frac{{2^p k}}{2}(1- \ch(2t/2^{n+p})w(1)^2
- \frac{2^p k}{2}\frac{ \left(\ch(2t/2^{n+p})-1 \right)^2 }{\ch(2t/2^{n+p})} 
\sum_{l =1}^{{2^p k -1}} 
  w\left(\frac{l}{{2^p k}}\right)^2 \right)
\\ \\
\displaystyle
 \exp  \left(-\frac{t}{2^{n+p}} \sum_{l =0}^{{2^p k-1}}
c\big(\frac{{2^p k}-l}{2^{n+p}}t ,
\sqrt{2^p k\ \sh(2t/2^{n+p})  } w\big(\frac{l}{{2^p k}}\big)
 + \frac{x}{  \ch(2t/2^{n+p})^{l } }\big) \right) 
\  dm_W(w).
\end{array}
\end{equation}
 The integrated terms are similar to the $F_n, G_n, H_n$
treated in Lemma \ref{convergences}. We get the limit of
  $v^{(t)}_{2^{n+p}}(\tau,x)$ by letting $p$ go to infinity in the integral and
 do not need the additional hypotheses. The only difference is that
$$
\lim_{p\rightarrow \infty} 2^pk \ \sh(2t/2^{n+p}) =  \tau.
$$
To find the limit of  $\partial_x^2 v^{(t)}_{2^n}(\tau,x)$,
one has to derivate (\ref{vn-tau-dyadique}) twice with respect to $x$.
The derivatives of the first exponential term contain the expression
$$
M= - \frac{\sh(2t/2^{n+p})}{\ch(2t/2^{n+p})}
 \sum_{l =0}^{{2^p k -1}}
 \left(  \sqrt{{2^p k}\ \sh(2t/2^{n+p})  }   w\left(\frac{l}{{2^p k}}\right)  
 + \frac{x}{\ch(2t/2^{n+p})^{l}} \right)
 \frac{1}{\ch(2t/2^{n+p})^{l}},
$$
as well as its square and its derivative (with respect to $x$).
We estimate $M$ by  $K( ||w||_{\infty}+|x|)$ with a constant $K$ depending only
on $t$. The other exponential terms and their derivatives are bounded with
respect to $w$. Therefore we can apply the dominated convergence Theorem,
since $\int_{C_W } ||w||^2 \ dm_W(w)$ is bounded (by Fernique's Theorem, see Section \ref{rappels} ).
As for the limits themselves, the same techniques can be applied as in
the proof of Lemma \ref{convergences}. Note that the estimate of 
$\partial^{3}_x c$
 is needed to treat terms containing  $\partial^2_x c$ for we use the mean
value theorem. \hfill $\square$

\vskip 0.2 cm

\begin{lem}\label{bornes} 
There exists $C>0$, depending on  $t$ and $x$, such that
$$
\begin{array}{ccc}
\displaystyle
\forall j\in \{0,\dots, 4\},\ 
\forall n\in \N^*,\ \forall k \in \{0,\dots, n -1\},\ 
\forall \tau \in 
 \left[\frac{kt}{n}, \frac{(k+1)t}{n}\right],\quad \\
\\\displaystyle
\left| \frac{\partial^j v^{(t)}_{n}(\tau,x)}{\partial x^j}\right| \leq C.
\end{array}
$$
\end{lem}

\noindent{\it Proof  of Lemma \ref{bornes}}\\
We first treat  the interval $ [\frac{2kt}{2n},\frac{(2k+1)t}{2n}]$.
Formula (\ref{vn-intermediaire}) shows that
$$
\label{vn(tau)}
\begin{array}{lllll}
\displaystyle
v_n(\tau,x)= \int_{\R^{k+1}}
(2\pi \sh(4(\tau-kt/n)))^{-1/2} (2\pi \sh(2t/n))^{-k/2}\\ \\
\displaystyle
 \exp\left(-\frac{1}{2}
\frac{\sh(4(\tau - kt/n ))}{\ch(4(\tau - kt/n ))}x^2 \right)
\exp\left(-\frac{1}{2}
\frac{\ch(4(\tau -kt/n))}{\sh(4(\tau -kt/n))}\sigma_1^2 \right)
\\ \\
\displaystyle
 \exp\left(-\frac{1}{2}\frac{\ch(2t/n)}{\sh(2t/n)} 
\sum_{j=1}^k 
\left( \sigma_{k+2-j} -\frac{ \sigma_{k+1-j} }{\ch(2t/n)} \right)^2  \right)
\\ \\
\displaystyle
\exp\left(-\frac{1}{2}\frac{\sh(2t/n)}{\ch(2t/n)}
\sum_{j=1}^k \left(\sigma_{k+1-j}+
\frac{x }{\ch(4(\tau-kt/n))\ch(2t/n)^{k-j}}\right)^2 \right)
 \\
\\
\displaystyle
\exp   \left(-\frac{t}{n}\sum_{j =1}^{k} c\left(\frac{j}{n}t ,
\sigma_{k+1-j} 
 + \frac{x}{ \ch(4(\tau -kt/n   ))
 \ch(2t/n)^{k-j} }\right)
 \right)
\\  \\
\displaystyle 
v_0\left(
\sigma_{k+1}
 + \frac{x}{ \ch(4(\tau - kt/n )) \ch(2t/n)^{k} }  \right) 
\quad d\sigma_1\dots d\sigma_{k+1}.
\end{array}
$$
The space-derivatives of order at most $4$ of  all terms but one are 
bounded by constants depending on  $t$, $||v_0^{(j)}||$
and $||\partial_x^jc||$ ($j\leq 4$), but not on $k$, $n$, $\tau$ or $x$.
The derivative of 
$$
 \exp\left(-\frac{1}{2}
\frac{\sh}{\ch}(4(\tau -\frac{ kt}{n} ))x^2 \right)
\exp\left(-\frac{1}{2}\frac{\sh}{\ch}(\frac{2t}{n})
\sum_{j=1}^k \left(\sigma_{k+1-j}+
\frac{x }{\ch(4(\tau-\frac{ kt}{n}))\ch(\frac{2t}{n})^{k-j}}
\right)^2 \right)
$$
is less easy to treat. It is the product of the exponential term and of 
$$
 -\frac{\sh}{\ch}(4(\tau -\frac{ kt}{n} ))x\
- \ 
\frac{\sh}{\ch}(\frac{ 2t}{n})
\sum_{j=1}^k \left(\sigma_{k+1-j}+
\frac{x }{\ch(4(\tau-\frac{ kt}{n}))\ch(\frac{2t}{n})^{k-j}}\right)
 \frac{1}{\ch(4(\tau-\frac{ kt}{n}))\ch(\frac{2t}{n})^{k-j}}.
$$
We then have to estimate expressions such as
$$
\left( \sum  a_sb_sc_s \right)^j \exp\left(- \frac{1}{2} \sum a_s b_s^2 \right),
$$
where the $a_s$ are quotients $\sh/\ch$ and the $c_s$ are
the $ (\ch(4(\tau-\frac{ kt}{n}))^{-1}\ch(\frac{2t}{n})^{-s}.$
Applying Cauchy-Schwartz inequality to the sum outside of the exponential term
shows it is smaller than
$ \sqrt{ \sum  a_sb_s^2  \sum  a_s c_s^2}  $. The first factor is absorbed by
 the exponential and the second one is a geometric sum, which is bounded.\\
It is eventually possible to estimate the derivatives of order
 $0\leq j\leq 4$ by
$$
\begin{array}{lllll}
\displaystyle
C \int_{\R^{k+1}}
(2\pi \sh(4(\tau-kt/n)))^{-1/2} (2\pi \sh(2t/n))^{-k/2}
\exp\left(-\frac{1}{2}
\frac{\ch(4(\tau -kt/n))}{\sh(4(\tau -kt/n))}\sigma_1^2 \right)
\\ \\
\displaystyle \qquad   \qquad
 \exp\left(-\frac{1}{2}\frac{\ch(2t/n)}{\sh(2t/n)} 
\sum_{j=1}^k 
\left( \sigma_{k+2-j} -\frac{ \sigma_{k+1-j} }{\ch(2t/n)} \right)^2  \right)
 \ d\sigma_1 \dots d\sigma_{k+1}\ ,
\end{array}
$$
where the constant $C$ depends on $t$, $||v_0^{(j)}||_{\infty}$ et
 $||\partial^j c||_{\infty}$ ($0\leq j\leq 4$) only. The integral is equal to
$$
\ch\left(4(\tau-\frac{ kt}{n})\right)^{-1/2}\ch(2t/n)^{-k/2}\leq 1,
$$
which establishes our claim for the interval 
$[\frac{2k}{2n}t, \frac{2k+1}{2n}]$.
\\
On the following interval,
 $[\frac{2k+1}{2n}t, \frac{2k+2}{2n}]$,
$$
v_n^{(t)}(\tau,x) =
\exp\left(-2(\tau-\frac{2k+1}{2n}t)
c(\frac{k+1}{n}t) , x) \right) v_n^{(t)}(\frac{2k+1}{2n}t,x) .
$$
It is the product of two functions having bounded space-derivatives (of order
 at most $4$). This proves the estimations on $[0,t]$. \hfill $\square$

\vskip 0.2 cm

\begin{lem}\label{cvu}
Let  $(u_n)_{n\in \N}$ be a sequence of functions continuous on $[0,T]$,
piecewise $C^1$ on $[0,T]$ and satisfying
$$
\exists C > 0 \ : \ 
\forall n\in \N,\ \forall \tau \in [0,T] \ {such \ that }\ 
u'_n(\tau) \ { exists},\quad
|u'_n(\tau)|\leq C.
$$
Suppose $D$ is a dense subset of $[0,T]$ and $u$, a continuous function on 
$[0,T]$, such that
$$
\forall \tau \in D, \ \lim_n u_n(\tau)= u(\tau).
$$
Then the sequence $(u_n)$ has a subsequence which converges to $u$ uniformly on
$[0,T]$.
\end{lem}
\vskip 0.2 cm
\noindent{\it Proof.}
The bounds on the derivatives show that  $(u_n)_{n\in \N}$ is bounded and 
equicontinuous on $[0,T]$. Ascoli's Theorem yields the existence of 
a uniformly converging subsequence. Its limit $\tilde{u}$ is continuous on
$[0,T]$ and equal to $u$ on the dense subset $D$, which justifies the Lemma.
 \hfill $\square$
\vskip 0.3 cm 
\noindent{\it Proof  of  Theorem \ref{solution-weak}.}\\
Let $x$ be fixed. We have proved that the (numerical) sequences
 $(v^{(t)}_{2^n}(\tau,x))_{n\in \N}$
 and $\displaystyle( (\partial_x)^2 v^{(t)}_{2^n}(\tau,x))_{n\in \N}$ converge
(respectively) to $v(\tau,x)$  and   $\partial_x^2v(\tau,x)$ 
for all $\tau$ belonging to the dense subset $\ddd$ defined in Lemma  
\ref{cv-sur-D}. Now to  prove the uniform convergence on $[0,t]$
we apply Lemma  \ref{cvu}.
The continuity of $v_{2^n}^{(t)}(\cdot,x)$ is a consequence of its definition.
The continuity of $\partial_x^2 v_{2^n}^{(t)}(\cdot,x)$ can be proved by 
studying the expressions appearing in the proof of Lemma \ref{bornes} 
and so is the derivability with respect to $\tau$.

To establish the bounds on the $\tau$-derivatives let us recall that
\begin{enumerate} 
\item
On the even intervals $[ 2kt/2.2^n, (2k+1)t/ 2.2^n]$, 
$$
\frac{\partial v_{2^n}^{(t)}(\tau,x)}{\partial \tau } =
2\frac{\partial^2 v_{2^n}^{(t)}(\tau,x)}{\partial x^2}-2x^2 v_{2^n}^{(t)}(\tau,x),
$$
\item
on the odd intervals $[ (2k+1)t/2.2^n, (2k+2)t/ 2.2^n]$, 
$$
\frac{\partial v_{2^n}^{(t)}(\tau,x)}{\partial \tau } 
= -2c((2k+2)t/ 2.2^n, x) v_{2^n}^{(t)}(\tau,x).
$$
\end{enumerate}
The bounds on the time derivatives come from the estimates on the space
 derivatives~: to show that $\partial_{\tau} v_{2^n}^{(t)} $ is bounded one
needs the $x$-derivatives up to order $2$, to treat 
$\partial_{\tau}(\partial_x)^2 v_{2^n}^{(t)}(\tau,x)$ one needs the $x$-
 derivatives up to order $4$. The bounds concerning the space-derivatives 
have been established in Lemma  \ref{bornes}. It follows that a subsequence of
 $( v_{2^n}^{(t)})$  (resp. of $(\partial_{\tau} v_{2^n}^{(t)})$ ) converges
uniformly on  $[0,t]$. Let us denote its indexes by $\ph(n)$.\\
Next we write more concisely the system of equations defining
$v_n^{(t)}$ (according to the interval).
On  $[0,t]\setminus \{ kt/2^n\}, 0\leq k\leq 2^n$, $v_n^{(t)}$ satisfies
$$
\frac{\partial v^{(t)}_{n}(\tau,x)}{\partial \tau } 
=
2\beta_n(\tau) \left( \frac{\partial^2 v^{(t)}_{n}(\tau,x)}{\partial x^2 } 
- x^2  v^{(t)}_{n}(\tau,x)\right) 
-2( 1-\beta_n(\tau))c_n(\tau,x) v^{(t)}_n(\tau,x),
$$
with
$$
\begin{array}{lllllll}
\beta_{n}(\tau) & = & 1 & {\rm on } & 
\displaystyle \left] \frac{2k}{2n}t,  \frac{2k+1}{2n}t \right]  
\\ \\
 & = & 0 & {\rm on } & 
\displaystyle \left] \frac{2k+1}{2n}t,  \frac{2k+2}{2n}t \right]    \\ \\
{\rm and }\ 
c_{n}(\tau,x) & = & c( \frac{2k+2}{2n}t,x)  & {\rm on } & 
\displaystyle \left] \frac{2k}{2n}t,  \frac{2k+2}{2n}t \right]  .  \\
\end{array}
$$
This equation still holds for the subsequence indexed by $\ph(n)$. 
The uniform convergence allows us to integrate on any subinterval $[0,s]$
of $[0,t]$ :
$$
\begin{array}{ll}
\displaystyle
v^{(t)}_{\ph(n)}(s,x) -v^{(t)}_{\ph(n)}(0,x) 
=\\
\displaystyle \quad
2 \int_0^{s}
\beta_{\ph(n)}(\tau)
 \left( \frac{\partial^2 v^{(t)}_{\ph(n)}(\tau,x)}{\partial x^2 } 
- x^2  v^{(t)}_{\ph(n)}(\tau,x)\right) 
-( 1-\beta_{\ph(n)}(\tau))c_{\ph(n)}(\tau,x) v^{(t)}_{\ph(n)}(\tau,x) \ d\tau.
\end{array}
$$
Now let $n$ tend to infinity. The following result (\cite{Kahane}) shows what
 become of $\beta_n$ and of the factor $2$ :
\begin{lem}\label{Lemme-Kahane}
Let  $\beta_n$ be the function  defined above. Let  $(\psi_n)$ be a sequence
 of functions belonging to $\in L^1([0,t])$ and suppose it converges uniformly
on $[0,t]$ to a limit $\psi$. Then, for all 
$0\leq \tau\leq \sigma\leq t$,
$$
\lim_{n\rightarrow \infty} \int_{\tau}^{\sigma}
\beta_n\psi_n \ ds = 
\lim_{n\rightarrow \infty} \int_{\tau}^{\sigma}
(1- \beta_n)\psi_n \ ds = \frac{1}{2} \int_{\tau}^{\sigma} \psi \ ds.
$$
\end{lem}

Eventually, for all $s\leq t$, we obtain
$$
v(s,x) -v_0(x) 
=
\int_0^{s}
\left( \frac{\partial^2 v(\tau,x)}{\partial x^2 } 
- x^2  v\tau,x)\right) 
-c(\tau,x) v(\tau,x)\ d\tau\ .
$$
As $v(s,x)$ does not depend on $t$ (the intermediates  $v_n^{(t)}$ 
depend on $t$ but not the limit), the function $v$ is a solution of (\ref{pb}).
\hfill $\square$


\section{Another proof of Feynman  Kac formula}\label{another}

Before completing the proof of Theorem \ref{solution} we shall see that a
 small modification of the method developped in Section\ref{suite}
 yields
the following expression for $v$.  For  sufficiently small $t$,
$$
\begin{array}{lll}
v(t,x)= & \displaystyle \int_{C_W} v_0(w(2t) +x) \
 \exp\left(-\frac{1}{2} \int_0^{2t} \left(x + w(s)\right)^2 \ ds\right)\\ \\
 & \displaystyle \qquad \qquad \exp\left(- \frac{1}{2}\int_0^{2t}
c\left(t-s/2, w(s)+ x\right)\ ds \right)\ dm_W(w)\ .
\end{array}
$$
With $u(t,x)=v(t/2,x)$, $u$ satisfies 
$$
\frac{\partial u}{\partial t }(t,x)  -\frac{1}{2}
\frac{\partial^2 u}{\partial x^2 }(t,x)  +\frac{1}{2}(x^2  + c(t/2,x))u(t,x)
=0,
$$
which explains the differences with the usual expression.

 The first point is that
the demonstration does not use the It\^ o  integral at all. 
What is, perhaps, more
 significant is that both expressions of $v$ are not linked by a ``classical''
 change of variable on Wiener space.

\noindent{\it Proof.}
It essentially follows the same steps as in Section \ref{suite}.
Starting from Proposition {\ref{integrales-iterees}}, we consider the time
 sequence $(t_k)_{0\leq k\leq n}$ with $t_k= k \frac{2t}{n}$ instead of
$t_k=k/n$.
This yields 
$$
\begin{array}{llll}
\displaystyle 
 v^{(t)}_n(t,x)  =
\int_{C_W}  
\left(\frac{2t}{n \sh(2t/n)} \right)^{n/2}
 v_0\left(  w(2t) + \frac{x}{\ch(2t/n)^n}\right)
 \\
\\
 \displaystyle
\exp\left( -\frac{1}{2}\frac{\ch(2t/n)}{\sh(2t/n)}     \sum_{j=1}^{n}
\left[w( 2jt/n) -\frac{w( 2(j-1)t/n)}{\ch(2t/n)}   \right]^2
\right)
\\  \\
\displaystyle
\exp\left( \frac{n}{4t}    \sum_{j=1}^{n}
\left[w( 2jt/n) -w( 2(j-1)t/n)  \right]^2
\right)
\\  \\
\displaystyle
\exp\left( -\frac{1}{2} \frac{\sh(2t/n)}{\ch(2t/n)}  
             \sum_{j=0}^{n-1} \left(  w( 2jt/n) + 
  \frac{x}{\ch(2t/n)^{j}}  \right)^2   \right) 
\\ \\
\displaystyle
\exp\left( -\frac{t}{n}  \sum_{j=0}^{n-1}
 c\left(\frac{(n-j) t}{n}, w(2jt/n) +\frac{x}{\ch(2t/n)^{j}} \right) 
    \right)\ 
dm_W(w).
\end{array}
$$
Most of the terms are bounded and converge as
 in Section {\ref{suite}} or even 
more easily. It just remains to treat 
$$
 -\frac{1}{2}\frac{\ch(2t/n)}{\sh(2t/n)}     \sum_{j=1}^{n}
\left[w( 2jt/n) -\frac{w( 2(j-1)t/n)}{\ch(2t/n)}   \right]^2
+ \frac{n}{4t}    \sum_{j=1}^{n}
\left[w( 2jt/n) -w( 2(j-1)t/n)  \right]^2 .
$$
This expression splits into $A_n+ B_n$ where
$$
A_n=
\left( \frac{n}{4t}  -\frac{1}{2}\frac{\ch(2t/n)}{\sh(2t/n)}\right) w(2t)^2
-   \sum_{j=1}^{n-1}  w( 2jt/n)^2 \frac{ (\ch(2t/n)-1)^2 }{2\ch(2t/n)\sh(2t/n)}
$$
and
$$
B_n = \left( \frac{n}{2t}  - \frac{1}{\sh(2t/n)}\right)
  \sum_{j=1}^{n-1}  w( 2jt/n)\left( w( 2jt/n)- w( 2(j-1)t/n)   \right).
$$
The first term is negative and converges to $0$. The second one can be
 estimated as follows
$$
\begin{array}{llll}
|B_n| &\leq & \displaystyle
\left| \frac{n}{2t}  - \frac{1}{\sh(2t/n)}\right|\frac{n}{2t}
\sqrt{\frac{2t}{n} \sum_{j=1}^{n-1}  w( 2jt/n)^2 }
\sqrt{\frac{2t}{n} \sum_{j=1}^{n-1} \left(w(2jt/n)-w(2(j-1)t/n)\right)^2  }\\ \\
 &\leq & \displaystyle
M \sqrt{ 2t}||w|| 
\sqrt{\frac{2t}{n} \sum_{j=1}^{n-1} \left(w(2jt/n)-w(2(j-1)t/n)\right)^2  }.\\ \\
\end{array}
$$
As $ \sum_{j=1}^{n-1} \left(w(2jt/n)-w(2(j-1)t/n)\right)^2$ is the quadratic
 variation of a Brownian motion, the subsequence for $n=2^p$    converges
 to $\sqrt{2t}$  and it is smaller than $4n||w||^2$.
To sum up, $|B_n|\leq 4tM||w||^2 $ and converges to $0$. This estimation 
and Fernique's theorem allow us  to use Lebesgue dominated convergence theorem,
 provided $t$ is small enough. \hfill $\square$


\section{Proof of Theorem \ref{solution}}\label{affaiblissement}

To get the optimal form of the Theorem, it remains to prove that formula
(\ref{def-de-v}) gives a solution of Problem (\ref{pb}) even if
$v_0$ and $c$ satisfy much weaker assumptions. This will be done by 
approximating general $v_0$ and $c$ by regular functions
 and  showing that 
the solution of the approximating problem converges to that of the 
real  problem.


\begin{prop}
For $v_0\in L^2(\R)$ and $c$ measurable and inferiorly bounded on 
$]0,\infty)\times \R$ we define, following formula (\ref{def-de-v}),
$$
\begin{array}{lll}
S(v_0,c)(t,x)=  & \displaystyle \int_{C_W} v_0(\sqrt{2t }\ w(1) +x) \
 \exp\left(-t\int_0^1 \left(x +\sqrt{2t}\ w(s)\right)^2 \ ds\right)\\ \\
 & \displaystyle \qquad \qquad \exp\left(-t\int_0^1
c\left(t(1-s),\sqrt{2t}\ w(s)+ x\right)\ ds \right)\ dm_W(w)\ .
\end{array}
$$
For all $0\leq \alpha < \beta <\infty$, $S(v_0,c)$ belongs to
$L^2([\alpha,\beta]\times \R)$. Moreover,
\begin{equation}\label{ineg-S}
\int_{[\alpha,\beta]\times \R  }|S(v_0,c)(t,x)|^2 \ dt dx \leq ||v_0||^2
\int_{\alpha}^{\beta} e^{-2t\inf(c)}\ dt\ .
\end{equation}

\end{prop}
\noindent{\it Proof .}\\
Let us consider
$$
\begin{array}{lll}
I & \displaystyle = \int_{\alpha}^{\beta} \int_{\R} \int_{C_W} 
| v_0(\sqrt{2t }\ w(1) +x)|^2 
\exp\left(-2t\int_0^1 \left(x +\sqrt{2t}\ w(s)\right)^2 \ ds\right)\\ \\
&  \displaystyle \ \times 
\exp\left(-2t\int_0^1
c\left(t(1-s),\sqrt{2t}\ w(s)+ x\right)\ ds \right)\ dm_W(w) \ dx dt \quad .
\end{array} 
$$
The first exponential factor is smaller than $1$ and the second one, than 
$\exp(-2t \inf(c))$. By Fubini's Theorem,
$$
\begin{array}{lll}
I 
 & \displaystyle \leq 
  \int_{\alpha}^{\beta} e^{-2t\inf(c)}  \int_{C_W}\left(  \int_{\R}
| v_0(\sqrt{2t }\ w(1) +x)|^2  \ dx \right) \ dm_W(w) dt\\ \\
&  \displaystyle \leq 
  \int_{\alpha}^{\beta} e^{-2t\inf(c)}  \int_{C_W}\left(  \int_{\R}
| v_0(\xi)|^2  \ d\xi \right) \ dm_W(w) dt\ ,\\ \\
\end{array}
$$
with $\xi= \sqrt{2t }\ w(1) +x$. Eventually,
$$
I\leq || v_0||^2  \int_{\alpha}^{\beta} e^{-2t\inf(c)} \ dt.
$$
This shows that the integral 

$$
\int_{C_W} 
| v_0(\sqrt{2t }\ w(1) +x)|^2 
e^{-2t\int_0^1 \left(x +\sqrt{2t}\ w(s)\right)^2 \ ds}
e^{-2t\int_0^1
c\left(t(1-s),\sqrt{2t}\ w(s)+ x\right)\ ds } \ dm_W(w) 
$$
converges for almost all $(t,x)\in [\alpha,\beta]\times \R$. By H\"older's
inequality it follows that $S(v_0,c)(t,x)$ is defined for
the same $(t,x)$ and that $S(v_0,c)$ satisfies the inequality (\ref{ineg-S}).
\hfill $\square $

\begin{prop}
Let
$(v_0^{(n)} )\in L^2(\R)^{\N}$ converge to $v_0$ in  $L^2(\R)$ and
 $(c^{(n)} ) \in L^2(\R^+\times \R)^{\N}  $ converge to $c$  $ \in
 L^2(\R^+\times \R)$.
Assume that the $c^{(n)}$ and $c$ have a common lower bound $\mu\in \R$.\\
Then, for all
$\alpha, \beta \in \R^+$ satisfying $0\leq \alpha \leq \beta <\infty$,
 $S(v_0^{(n)}, c^{(n)})$ converges to $S(v_0,c)$ in 
 $L^2([\alpha,\beta]\times \R)$.
\end{prop}
\noindent{\it Proof .}
Let us introduce the intermediate $S(v_0,c^{(n)})$. Then by the preceding
Proposition 
$$
|| S(v_0,c^{(n)})- S(v_0^{(n)},c^{(n)})||_{L^2([\alpha,\beta]\times \R)}
\leq || v_0-v_0^{(n)}||_{L^2(\R)}
 \left(\int_{\alpha}^{\beta} e^{-2t\mu} \ dt \right)^{1/2}
$$
and this term converges to $0$. \\
The second term 
$ S(v_0,c^{(n)})- S(v_0, c)$ is more delicate. Let
$$\begin{array}{lll}
I_n(t,x)= & \displaystyle 
 \int_{C_W} v_0(\sqrt{2t } w(1) +x) 
 e^{-t\int_0^1 \left(x +\sqrt{2t} w(s)\right)^2  ds} \\ \\
& \displaystyle  \left(
e^{-t\int_0^1c^{(n)}\left(t(1-s),\sqrt{2t} w(s)+ x\right)  ds}
- e^{-t\int_0^1c\left(t(1-s),\sqrt{2t} w(s)+ x\right)  ds}\right)
  dm_W \ ,
\end{array}
$$
so that we may write
$$
 ||   S(v_0,c^{(n)})- S(v_0, c) ||^2_{L^2([\alpha,\beta]\times \R)}=
\int_{\alpha}^{\beta} \int_{\R} I_n(t,x)^2 \ dt \ dx. 
$$
We shall prove that the sequence $(I_n(t,x))_n$ converges to $0$ when $n$ 
goes to infinity and that it is uniformly bounded by a function
$g(t,x)$ belonging to  $L^2([\alpha,\beta]\times \R)$. Then, thanks to the
 dominated
 convergence Theorem, 
$||   S(v_0,c^{(n)})- S(v_0, c) ||^2_{L^2([\alpha,\beta]\times \R)}$ will 
converge to $ 0$.\\
As usual, the first exponential term of $I_n(t,x)$ is smaller than $1$. 
Both 
$-t\int_0^1c^{(n)}(t(1-s),\sqrt{2t} w(s)+ x)  ds$
and 
$-t\int_0^1c(t(1-s),\sqrt{2t} w(s)+ x)  ds$
being smaller than $-t\mu$, 
$$
\begin{array}{ll}
\displaystyle   \left|
e^{-t\int_0^1c^{(n)}\left(t(1-s),\sqrt{2t} w(s)+ x\right)  ds}
- e^{-t\int_0^1c\left(t(1-s),\sqrt{2t} w(s)+ x\right)  ds}\right|\\ \\
 \displaystyle   \leq 
\left|t\int_0^1(c^{(n)}-c)(t(1-s),\sqrt{2t}\  w(s)+ x)ds\right|e^{-t\mu}
\end{array}
$$
and we can estimate $|I_n|$ by
$$
\begin{array}{ccc}
|I_n(t,x)|\leq t e^{-t\mu}
 &\displaystyle \underbrace{ \int_{C_W}| v_0(\sqrt{2t } w(1) +x)| 
\int_0^1\left|(c^{(n)}-c)(t(1-s),\sqrt{2t}\  w(s)+ x)\right| ds   dm_W}\ .\\
& M_n(t,x)\\
\end{array}
$$
Thanks to (\ref{Wiener-nbfinipoints}), $M_n(t,x)$ can be written as an 
integral on $\R$ :
$$
\begin{array}{ll}
\displaystyle  
M_n(t,x) & = \displaystyle  
 \int_0^1 \int_{\R^2}| v_0(\sqrt{2t }\xi_2 +x)|
\left|(c^{(n)}-c)(t(1-s),\sqrt{2t}\  \xi_1+ x)\right| 
f_{s,1}(\xi_1, \xi_2 )  d\xi_1d\xi_2  \ ds ,\\ \\
\end{array}
$$
where $f_{s,1}$ is the gaussian density. The change of variables
$$
u= \sqrt{2t }\xi_2 +x,\ v=t(1-s),\ w=\sqrt{2t}\  \xi_1+ x
$$
gives
$$
\begin{array}{ll}
\displaystyle  
M_n(t,x) & = \displaystyle  
\frac{1}{4\pi t} \int_{\R^3} \frac{{\mathbf 1}_{[0,t]}(v)}{\sqrt{v(t-v)}}   
\left|(c^{(n)}-c)(v,w)\right| e^{  -\frac{(w-x)^2}{4(t-v)}}
 e^{-\frac{(u-w)^2}{4v}}|v_0(u)|\ du dv dw.
\end{array}
$$
As 
$$
 \int_{\R} e^{-\frac{(u-w)^2}{4v}}|v_0(u)|\ du \leq ||v_0||_{L^2(\R)}
 (2\pi v)^{1/4},
$$
$M_n(t,x)$ is smaller than
$$
 \frac{1}{4\pi t}||v_0||_{L^2(\R)}(2\pi)^{1/4}
||c-c^{(n)}||_{L^2(\R^+\times \R)}
 \sqrt{ \int_{\R^2}   \frac{  {\mathbf 1}_{[0,t]}(v)  }{v(t-v)}    v^{1/2}  
 \exp\left(  -\frac{(w-x)^2}{2(t-v)}\right) \ dv dw }.
$$
The integral appearing in the square root converges and does not depend on 
$n$, which shows that $M_n(t,x) $ and $I_n(t,x)$ go to 
$0$.\\
Now for the uniform estimate. Clearly
$$
|I_n(t,x)|\leq 
 \int_{C_W}| v_0(\sqrt{2t } w(1) +x) | (e^{-t\mu}+ e^{-t\mu})  dm_W := g(t,x).
$$
The function $g$ is in  $L^2([\alpha,\beta]\times \R)$  since
$$
\begin{array}{lll}
\displaystyle \int_{\alpha}^{\beta}\int_{\R} g(t,x)^2 \ dx dt &
 \displaystyle \leq  \int_{\alpha}^{\beta}
4e^{-2t\mu}\int_{\R}
\left(\int_{C_W}|v_0(\sqrt{2t}w(1)+x)|dm_W\right)^2\ dx dt\\ \\
&\displaystyle
\leq  \int_{\alpha}^{\beta}
4e^{-2t\mu}\int_{\R}\int_{C_W}|v_0(\sqrt{2t}w(1)+x)|^2dm_W\ dxdt\\ \\
&\displaystyle
\leq  \int_{\alpha}^{\beta}
4e^{-2t\mu}
\int_{C_W} \int_{\R}|v_0(\sqrt{2t}w(1)+x)|^2 \ dx\  dm_W dt \ . \\ \\
\end{array}
$$
The change of variables $\xi =\sqrt{2t}w(1)+x$ allows to write
$$
\int_{\alpha}^{\beta}\int_{\R} g(t,x)^2 \ dx dt \leq ||v_0||^2_{L^2(\R)}
\leq  \int_{\alpha}^{\beta} 4e^{-2t\mu} dt < \infty ,
$$
which concludes the proof. \hfill $\square$\\

Any $v_0\in L^2(\R)$ can be approximated (in $L^2(\R)$) by a sequence
 $(v_0^{(n)}) $ of
functions satisfying the hypotheses of Theorem \ref{solution-weak}. Similarly,
any function $c$ belonging to
 $L^2(]0,\infty) \times \R)\cap L^{\infty}(]0,\infty)\times \R)$
 is the limit (in 
 $L^2(]0,\infty)\times \R)$) of a sequence  $(c^{(n)})_n $ of functions
satisfying  the hypotheses of Theorem \ref{solution-weak}. Moreover we can
suppose
this sequence to be bounded in $ L^{\infty}(]0,\infty)\times \R)$.
Then, according to Theorem \ref{solution-weak},
$v_n :=S(v_0^{(n)}, c^{(n)})$ is a solution of
$$
\left\{
\begin{array}{lll}
\displaystyle
\frac{\partial v_n}{\partial t }(t,x)  -
\frac{\partial^2 v_n}{\partial x^2 }(t,x)  +(x^2  + c^{(n)}(t,x))v_n(t,x)
=0 \quad {\rm on }\quad ]0,\infty[\times \R 
 \\
\\ \displaystyle
v_n(0,x)=v_0^{(n)}(x) \ .
\end{array}\right.
$$
Let $\ph $ be a smooth test function on $]0,\infty[\times \R$. Integrations
 by part give
$$
<v_n,\frac{\partial\ph}{\partial t} > -
<v_n,\frac{\partial^2 \ph}{\partial x^2} >
+ <x^2 v_n,\ph> +< c^{(n)} v_n,\ph>=0,
$$
where the brackets stand for $L^2(]0,\infty[\times \R)$ products.
As $v_n$ converges to $S(v_0,c)$ in any $L^2([\alpha,\beta]\times \R)$, 
we obtain
$$
<S(v_0,c),\frac{\partial\ph}{\partial t} > -
<S(v_0,c),\frac{\partial^2 \ph}{\partial x^2} >
+ <x^2 S(v_0,c),\ph> +< c S(v_0,c),\ph>=0.
$$
Moreover, 
$v_n(0,\cdot) = v^{(n)}$ converges to $v_0$, which gives the equality
$S(v_0,c)(0,\cdot)= v_0$.
Hence $S(v_0,c)$ is a solution of (\ref{pb}) in a weak sense.

\appendix

\section{Computation of some Wiener integrals}\label{explicit-expressions}

Another consequence of Theorem \ref{solution} is the following Proposition :
\begin{prop}\label{non-perturbee}
Suppose  $v_0$  satisfies, for all positive $t$ and real $x$ 
$$
\int_{\R}
(2\pi\ \sh(2t))^{-1/2} 
\exp\left(-\frac{1}{2}\frac{\ch(2t)}{\sh(2t)} y^2 
 -\frac{1}{2} \frac{\sh(2t)}{\ch(2t)} x^2 \right)
\left| v_0\left( y +\frac{x}{\ch(2t)} \right) \right| \ dy
<\infty.
$$
Then 
for all  $t>0$ and  $x\in \R$ we can write
$$
\begin{array}{lll}
\displaystyle \int_{C_W} v_0(\sqrt{2t}\ w(1)+x)\ \exp\left(-t 
\int_0^1 \left(x+\sqrt{2t}\ w(s) \right)^2\ ds\right) \ dm_W(w) \\ \\
=\qquad \displaystyle  \int_{\R}
(2\pi\ \sh(2t))^{-1/2} 
\exp\left(-\frac{1}{2}\frac{\ch(2t)}{\sh(2t)} y^2 
 -\frac{1}{2} \frac{\sh(2t)}{\ch(2t)} x^2 \right)
\left| v_0\left( y +\frac{x}{\ch(2t)} \right) \right| \ dy \ .
\end{array}
$$
\end{prop}
\vskip 0.3 cm 
\noindent{\it Proof.} 
Suppose $v_0$ belongs to  $L^2(\R)$. Then Theorem \ref{solution} shows 
that, as the perturbation $c$ is equal to $0$, the left hand side is 
the solution of the heat equation for the harmonic oscillator, with initial 
condition $v_0$.  The right hand side is the solution  $U_tv_0$ 
of the same problem,
given by Mehler's formula. The equality follows.\\
When $v_0$ is not in $L^2(\R)$, the equality holds
 for $v_0\ph_n$ where $\ph_n$ is  a convenient truncature. 
Then the assumption on $v_0$ allows to use the Theorems of dominated and of
 monotone convergence.  \hfill $\square$
\vskip 0.3 cm

When  $v_0(x)=1, x$ or $x^2$, it is easy to compute $U_tv_0$ and to deduce
the following equalities from 
these computations :
\begin{coro}\label{expr-k}
For $v_0 = 1$ we obtain
$$
k(t,x):= \int_{C_W}
e^{-t\int_0^1\left(x+\sqrt{2t}\ w(s)\right)^2 \ ds}\ dm_W(w)
=\frac{1}{\sqrt{\ch(2t)}}
\exp\left(- \frac{1}{2}\frac{\sh(2t)}{\ch(2t)}x^2\right).
$$
The case $v_0(x)=x$ gives 
$$
\int_{C_W}  w(1)  \exp\left( -t\int_0^1 \left(  x + \sqrt{2t}\ w(s)\right)^2 
\ ds\right) \ dm_W(w)
=\frac{x(1-\ch(2t))}{\sqrt{2t}\  \ch(2t)} \ k(t,x) 
$$
and for $v_0(x)=x^2$ we get
$$
\int_{C_W}  w(1)^2  \exp\big( -t\int_0^1 \big(  x + \sqrt{2t}\ w(s)\big)^2 
\ ds\big) \ dm_W(w) = 
\frac{1}{ 2t} 
\left( \frac{(1-\ch(2t))^2 x^2}{\ch^2(2t)} + \frac{\sh(2t)}{\ch(2t)} \right)\
k(t,x) \ .
$$
\end{coro}

As the perturbation $c$ is equal to $0$, \cite{Donsker-Lions} and 
\cite{Kuo} prove that, under certain conditions, the l.h.s. is the solution
of the heat equation problem for $H$. Nevertheless these integrals
are not mentioned explicitly in the literature, to the author's knowledge.


\end{document}